\documentclass[12pt, reqno]{amsart}
\usepackage{amssymb, mathrsfs}
\usepackage{latexsym}
\usepackage{amsmath}
\usepackage{amsthm}
\usepackage{amsfonts}
\usepackage{dsfont, upgreek}
\usepackage{mathtools}
\usepackage{epsfig}
\usepackage{amscd}
\usepackage{graphicx}
\usepackage{tikz-cd}
\usepackage{appendix}
\usepackage{enumerate}
\usepackage{color}
\usepackage{todonotes}

\usepackage[left=1.4in,right=1.4in,top=1.2in,bottom=1.2in]{geometry}

\usepackage{tikz}
\usetikzlibrary{decorations.markings}
\usetikzlibrary{arrows.meta}

\usepackage{extarrows}

\usepackage[colorlinks=true,linkcolor=blue,citecolor=blue]{hyperref}

\usepackage[all,cmtip]{xy}
\theoremstyle{plain}
\newtheorem{theorem}{Theorem}[section]

\newtheorem{claim}{Claim}

\theoremstyle{definition}

\theoremstyle{remark} 
\newtheorem{remark}[theorem]{Remark}

\numberwithin{equation}{section}

\newcommand{\ind}{\textup{Ind}}

\newcommand{\dist}{\mathrm{dist}}
\newcommand{\sgn}{\mathrm{sgn}}
\newcommand{\R}{\mathbb{R}}

\newcommand{\Sc}{\mathrm{Sc}}

\newcommand{\cl}{\mathrm{C}\ell}

\newcommand{\hotimes}{\mathbin{\widehat{\otimes}}}

\makeatletter
\let\save@mathaccent\mathaccent
\newcommand*\if@single[3]{%
	\setbox0\hbox{${\mathaccent"0362{#1}}^H$}%
	\setbox2\hbox{${\mathaccent"0362{\kern0pt#1}}^H$}%
	\ifdim\ht0=\ht2 #3\else #2\fi
}
\newcommand*\rel@kern[1]{\kern#1\dimexpr\macc@kerna}
\newcommand*\overbar[1]{\@ifnextchar^{{\wide@bar{#1}{0}}}{\wide@bar{#1}{1}}}
\newcommand*\wide@bar[2]{\if@single{#1}{\wide@bar@{#1}{#2}{1}}{\wide@bar@{#1}{#2}{2}}}
\newcommand*\wide@bar@[3]{%
	\begingroup
	\def\mathaccent##1##2{%
		\let\mathaccent\save@mathaccent
		\if#32 \let\macc@nucleus\first@char \fi
		\setbox\z@\hbox{$\macc@style{\macc@nucleus}_{}$}%
		\setbox\tw@\hbox{$\macc@style{\macc@nucleus}{}_{}$}%
		\dimen@\wd\tw@
		\advance\dimen@-\wd\z@
		\divide\dimen@ 3
		\@tempdima\wd\tw@
		\advance\@tempdima-\scriptspace
		\divide\@tempdima 10
		\advance\dimen@-\@tempdima
		\ifdim\dimen@>\z@ \dimen@0pt\fi
		\rel@kern{0.6}\kern-\dimen@
		\if#31
		\overline{\rel@kern{-0.6}\kern\dimen@\macc@nucleus\rel@kern{0.4}\kern\dimen@}%
		\advance\dimen@0.4\dimexpr\macc@kerna
		\let\final@kern#2%
		\ifdim\dimen@<\z@ \let\final@kern1\fi
		\if\final@kern1 \kern-\dimen@\fi
		\else
		\overline{\rel@kern{-0.6}\kern\dimen@#1}%
		\fi
	}%
	\macc@depth\@ne
	\let\math@bgroup\@empty \let\math@egroup\macc@set@skewchar
	\mathsurround\z@ \frozen@everymath{\mathgroup\macc@group\relax}%
	\macc@set@skewchar\relax
	\let\mathaccentV\macc@nested@a
	\if#31
	\macc@nested@a\relax111{#1}%
	\else
	\def\gobble@till@marker##1\endmarker{}%
	\futurelet\first@char\gobble@till@marker#1\endmarker
	\ifcat\noexpand\first@char A\else
	\def\first@char{}%
	\fi
	\macc@nested@a\relax111{\first@char}%
	\fi
	\endgroup
}
\makeatother

\makeatletter
\newcommand{\proofpart}[2]{%
	\par
	\addvspace{\medskipamount}%
	\noindent\emph{Part #1: #2}\par\nobreak
	\addvspace{\smallskipamount}%
	\@afterheading
}
\makeatother
\begin{document}
	
\title[A proof of Gromov's cube inequality]{ A proof of Gromov's cube inequality on scalar curvature }

\author{Jinmin Wang}
\address[Jinmin Wang]{School of Mathematical Sciences and Shanghai Center for Mathematical Sciences, Fudan University}
\email{wangjinmin@fudan.edu.cn}
\thanks{The first author is partially supported by NSFC11420101001.}
\author{Zhizhang Xie}
\address[Zhizhang Xie]{ Department of Mathematics, Texas A\&M University }
\email{xie@math.tamu.edu}
\thanks{The second author is partially supported by NSF 1800737 and 1952693.}
\author{Guoliang Yu}
\address[Guoliang Yu]{ Department of
	Mathematics, Texas A\&M University}
\email{guoliangyu@math.tamu.edu}
\thanks{The third author is partially supported by NSF 1700021, 2000082, 2247313, and the Simons Fellows Program.}

\date{}
\begin{abstract}
Gromov proved a cube inequality on the bound of distances between opposite faces of a cube equipped with a positive scalar curvature metric in dimension $\leq 8$ using  minimal surface method. He conjectured that the cube inequality also holds in dimension $\geq 9$. In this paper, we prove Gromov's cube inequality in all dimensions with the optimal constant via Dirac operator method. In fact, our proof yields a strengthened version of Gromov's cube inequality, which does not seem to be accessible by the minimal surface method. 
\end{abstract}
\maketitle

\section{Introduction}

In \cite{Gromov:2019aa}, Gromov proved a cube inequality on scalar curvature in dimension $\leq 8$ using minimal surface method \cite{RSSY79b, RSSY79, Schoen:2017aa}.
 In dimension $\geq 9$,  Gromov provided an approach toward a proof of this inequality  based on  unpublished 
results  of Lohkamp \cite{Lohkamp:2018wp} or a generalization of   Schoen-Yau's result  in \cite{Schoen:2017aa}. He stated that the cube inequality should be regarded as conjectural  in dimension $\geq 9$, cf.  \cite[Section 5.2, page 255]{Gromov:2019aa}.
 On the other hand, by applying a quantitative relative index theorem, the second author gave an index theoretic proof of Gromov's cube inequality  in all dimensions but with a suboptimal constant \cite{Xie:2021tm}. In this paper, we give an index theoretic proof of a strengthened Gromov's cube inequality with the optimal constant in all dimensions. 

\begin{theorem}\label{thm:cube}
	Let $g$ be a Riemannian metric on the cube $I^n=[0,1]^n$. Suppose $\ell_i=\dist(\partial_{i-},\partial_{i+})$ is the $g$-distance between the $i$-th pair of opposite faces $\partial_{i-}$ and $\partial_{i+}$ of the cube. Then the following hold. 
	\begin{enumerate}
		\item If the scalar curvature of $g$ satisfies $\Sc_g\geq k>0$, then 
		\begin{equation}\label{eq:cube}
		\sum_{i=1}^{n}\frac{1}{\ell_i^2}\geq \frac{kn}{4\pi^2(n-1)}.  
		\end{equation}
		 Consequently, we have 
		\begin{equation}\label{eq:cube2}
		\min_{1\leq i \leq n} \dist(\partial_{i-},\partial_{i+}) \leq 2\pi\sqrt{\frac{n-1}{k}}.
		\end{equation}
		\item If the scalar curvature of $g$ satisfies $\Sc_g\geq k>0$ and the dihedral angles of $(I^n,g)$ (viewed as a manifold with corners) are $< \pi$, then \begin{equation}\label{eq:cubestrict}
			\sum_{i=1}^{n}\frac{1}{\ell_i^2}>\frac{kn}{4\pi^2(n-1)}.  
		\end{equation}
		Consequently, we have 
		\begin{equation}\label{eq:cubestrict2}
		\min_{1\leq i \leq n} \dist(\partial_{i-},\partial_{i+}) < 2\pi\sqrt{\frac{n-1}{k}}.
		\end{equation}
	\end{enumerate}
\end{theorem}

We remark that the strict inequality from Part (2) of the above theorem  holds  when $(I^n, g)$ is  a manifold with smooth boundary. In this case, the condition on dihedral angles becomes vacuous, since there are no dihedral angles when the boundary is smooth. More generally, the strict inequality from Part (2) is also applicable when the dihedral angles of $(I^n,g)$ are $\leq \pi$ somewhere. See Remark \ref{rm:anglePi}.


The above inequality strengthens Gromov's cube inequality in \cite{Gromov:2019aa}, where his inequality sign is stated as  $\geq$. 
We remark that the constant $4\pi^2$ in the inequality \eqref{eq:cube} is optimal.
This can be seen by considering warped product metrics on large flat $(n-1)$-dimensional cubes over an interval (see Remark \ref{rk:optimal} for more details). To prove the strict inequality, we use the index theory for manifolds with corners developed in \cite{Wang:2021tq} together with some techniques developed in \cite{WangXie2022}.

  
The proof of the above theorem is based on an explicit index computation of  the Dirac operator with an appropriately chosen potential on $\R^n$.  Our proof also works for a more general version of Gromov's cube inequality in the spin case. More precisely, the following Gromov's $\square^{n-m}$ inequality holds. 

\begin{theorem}\label{thm:cubespin}
	Let $X$ be an $n$-dimensional  compact connected spin manifold with corners.
	Suppose 
$f\colon X\to [-1, 1]^m$
	is  a continuous map that preserves corner structures.\footnote{This means $f$ maps codimension $k$ faces to codimension $k$ faces. }  Let $\partial_{i\pm}, i = 1, \dots, m$,  be the pullbacks of the pairs of the opposite faces of the cube $[-1, 1]^m$.
	Suppose  $Y_{\pitchfork}$ is an $(n-m)$-dimensional closed submanifold \textup{(}without boundary\textup{)} in $X$ that satisfies the following conditions:
	\begin{enumerate}[$(1)$]
		\item $\pi_1(Y_\pitchfork) \to \pi_1(X)$ is injective;
		\item $Y_{\pitchfork}$  is the transversal intersection\footnote{In particular, this implies that the normal bundle of $Y_{\pitchfork}$ is trivial.}  of $m$ orientable hypersurfaces $\{Y_i\}_{i=1}^m $ of $X$, where  each $Y_i$  separates $\partial_{i-}$ from $\partial_{i+}$;
		\item the higher index $\ind_{\Gamma}(D_{Y_\pitchfork})$ does not vanish in $KO_{n-m}(C^\ast_{\max}(\Gamma; \mathbb R))$, where ${\Gamma = \pi_1(Y_\pitchfork)}$ and $C^\ast_{\max}(\Gamma; \mathbb R)$ is the  maximal group $C^\ast$-algebra of $\Gamma$ with real coefficients.  
	\end{enumerate} 
	Then the following hold. 
	\begin{enumerate}
		\item If $\Sc(X) \geq  k>0$, then  the distances $\ell_i = \dist(\partial_{i-}, \partial_{i+})$ satisfy the following inequality:
		\begin{equation}\label{eq:spin1}
	\sum_{i=1}^m \frac{1}{\ell_i^2} \geq \frac{kn}{4\pi^2(n-1)}. 
	\end{equation}
	Consequently, we have 
		 \begin{equation}\label{eq:spin2}
	\min_{1\leq i \leq m} \dist(\partial_{i-}, \partial_{i+}) \leq  
	2\pi\sqrt{\frac{m(n-1)}{kn}}.
	\end{equation}
	\item Suppose in addition $\pi_1(Y_\pitchfork) \to \pi_1(X)$ is split injective. If  $\Sc(X) \geq  k>0$ and the dihedral angles of $X$ are $ < \pi$, then \begin{equation}
		\sum_{i=1}^m \frac{1}{\ell_i^2} > \frac{kn}{4\pi^2(n-1)}. 
		\end{equation}
		Consequently, we have 
			 \begin{equation}
		\min_{1\leq i \leq m} \dist(\partial_{i-}, \partial_{i+}) < 
		2\pi\sqrt{\frac{m(n-1)}{kn}}.
		\end{equation}
\end{enumerate}
\end{theorem}
Again we remark that the strict inequality from Part (2) of the above theorem  is applicable to the special case where $X$  a Riemannian manifold with smooth boundary. Indeed, the condition on dihedral angles becomes vacuous in the case where the boundary is smooth.

We point out that Cecchini \cite{MR4181824} and Zeidler \cite{Zeidler:2019aa,MR4181525} proved a special case of Theorem \ref{thm:cubespin} when $m=1$, i.e.,  Gromov's bandwidth conjecture. The cube inequality plays an important  role in the  recent proof of nonexistence of positive scalar curvature metrics on $5$-dimensional aspherical manifolds in \cite{Gromov:2020aa} (cf. another proof of the same result in \cite{Chodosh:2020tk}).

 We would like to thank Xianzhe Dai and Guofang Wei for very helpful comments.

\section{Proof of Gromov's cube inequality}

In this section, we prove Theorem \ref{thm:cube}. Our proof is inspired in part by Zeidler's proof for Gromov's bandwidth conjecture \cite[Theorem 1.4]{MR4181525}.  
Throughout this section, we assume the constant $k$ in Theorem \ref{thm:cube} to be $n(n-1)$, while the general case follows by rescaling the metric.
\begin{proof}[Proof of Theorem \ref{thm:cube}]

To avoid distraction of the main ideas, we first prove the $\geq$ inequality.
We will give the proof of the strict inequality in Part 2. 

\proofpart{1}{$\geq$ inequality}

	Assume to the contrary that
\begin{equation}\label{eq:contrary}
\sum_{i=1}^{n}\frac{1}{\ell_i^2}< \frac{n^2}{4\pi^2}.
\end{equation}
Then there exists $\varepsilon>0$ such that
\begin{equation}\label{eq:epsilon}
\sum_{i=1}^{n}\frac{\pi^2(1+\varepsilon)^2}{(\ell_i-2\varepsilon)^2}< \frac{n^2}{4}.
\end{equation}

Now we shall extend the metric on the cube $I^n$ to a complete metric on $\R^n$.   
We first give a construction of the metric for $n=2$, and then briefly describe the general case.

\begin{figure}[h]
	\centering
	\begin{tikzpicture}
	\draw (0.5,0.5) -- (0.5,-0.5) -- (-0.5,-0.5) -- (-0.5,0.5) -- (0.5,0.5);
	\draw (0.5,0.5) .. controls (1,0.5) and (2,1.3) .. (4,1.3);
	\draw (0.5,-0.5) .. controls (1,-0.5) and (2,-1.3) .. (4,-1.3);
	\draw (-0.5,0.5) .. controls (-1,0.5) and (-2,1.3) .. (-4,1.3);
	\draw (-0.5,-0.5) .. controls (-1,-0.5) and (-2,-1.3) .. (-4,-1.3);
	\draw (0.5,0.5) .. controls (0.5,1) and (1.3,2) .. (1.3,3);
	\draw (0.5,-0.5) .. controls (0.5,-1) and (1.3,-2) .. (1.3,-3);
	\draw (-0.5,0.5) .. controls (-0.5,1) and (-1.3,2) .. (-1.3,3);
	\draw (-0.5,-0.5) .. controls (-0.5,-1) and (-1.3,-2) .. (-1.3,-3);
	\filldraw (4,0) node {$X_2$};
	\filldraw (0,3) node {$X_1$};
	\filldraw (0,0) node {$X$};
	\end{tikzpicture}
	\caption{Extend the metric from $[0,1]^2$ to $X_1\cup X_2$ and rescale the metric.}
	\label{fig:cross}
\end{figure}
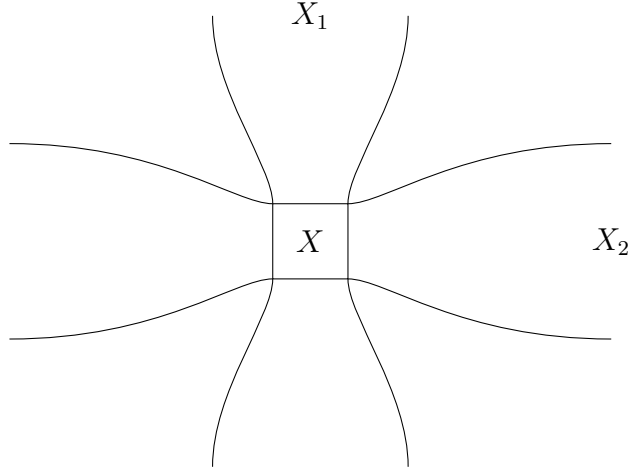

Denote $X=[0,1]^2$, $X_2=\R\times[0,1]$ and $X_1=[0,1]\times\R$.
Note that $X_1\cap X_2=X$, and $X_1\cup X_2$ has a shape of cross (see Figure \ref{fig:cross}).
 By the Whitney's extension theorem \cite{MR1501735} and a smooth partition of unity, there is a smooth metric $g_1$ on $X_1\cup X_2$ such that $g_1=g$ on $X$ and $g_1=dx_1^2+dx_2^2$ outside a compact set of $X_1\cup X_2$.

Suppose that $\Psi$ is a smooth positive function on $X_1\cup X_2$ such that $\Psi(x_1,x_2)$ is equal to $1$ if $(x_1,x_2)\in X$, and $\Psi(x_1,x_2)$ is a constant when $x_1^2+x_2^2$ is sufficiently large. Set $g_2=\Psi\cdot g_1$. By choosing $\Psi$ to be large enough outside $X$, we obtain a smooth metric $g_2$ on $X_1\cup X_2$ with the following properties (see Figure \ref{fig:cross}).

\begin{enumerate}[$(1)$]
	\item The $g_2$-distance $\dist_{g_2}(\partial_- X_1, \partial_+ X_1)$ between the two boundary components $\partial_- X_{1} = \{0\}\times \mathbb R$ and $\partial_+X_1 = \{1\}\times \R$ is $\geq\ell_1-\varepsilon/3$, where $\varepsilon$ is given in line \eqref{eq:epsilon}.
	\item The $g_2$-distance $\dist_{g_2}(\partial_- X_2, \partial_+ X_2)$ between the two boundary components $\partial_- X_{2} =  \mathbb R \times \{0\}$ and $\partial_+X_2 = \R\times \{1\}$ is $\geq \ell_2-\varepsilon/3$.
	\item $g_2=R^2(dx_1^2+dx_2^2)$ on $X_1\cup X_2$ outside a compact set for some $R>0$.
\end{enumerate}

Finally, using Whitney's extension theorem, we extend the metric $g_2$ on $X_1\cup X_2$ to obtain a complete metric $\overbar g$ on $\R^2$ such that $\overbar g=R^2(dx_1^2+dx_2^2)$ on $\mathbb R^2$ outside a compact subset. Since $\overbar g$ is flat outside a compact set, there exists $\sigma>0$ such that $\Sc_{\overbar g}\geq -\sigma$ on $\mathbb R^2$.

For $x=(x_1,x_2)\in\R^2$, we define
$$\varphi_1(x)=\sgn(x_1)\cdot \dist_{\overbar g}(x,\{0\}\times\R)\textup{  and }
\varphi_2(x)=\sgn(x_2)\cdot\dist_{\overbar g}(x,\R\times\{0\}),$$
where $\sgn$ is the sign function. 
Set $\psi_i=\varphi_i-\ell_i/2$. For each $i=1,2$, let $z_i$ be a smooth approximation of $\psi_i$ such that  $|z_i-\psi_i|\leq\varepsilon/3$ and
$\|\nabla z_i\|\leq 1+\varepsilon$  (cf. \cite[Proposition 2.1]{MR532376}). Here $\nabla$ is the gradient with respect to the metric $\overbar g$.  

To summarize, we have constructed a complete metric $\overbar g$  and smooth functions $z_1$ and $z_2$ on $\R^2$ with the following properties.
\begin{itemize}
	\item The metric $\overbar g$ restricts to the metric $g$ on $[0,1]^2$, and coincides with the Euclidean metric on $\R^2$ outside a compact set. In particular, there exists a positive number $\sigma>0$ such that $\Sc_{\overbar g}\geq -\sigma$ on $\R^2$.
	\item For each $x\in\R^2$, if both $|z_1(x)|\leq \ell_1/2-\varepsilon$ and $|z_2(x)|\leq  \ell_2/2-\varepsilon$, then $x$ lies in the cube $[0,1]^2$.
	\item For each $C > 0$, the subset $\{x\in \R^2: z_1(x)^2+z_2(x)^2\leq C\}$ is compact.
	\item $\|\nabla z_i\|\leq 1+\varepsilon$, where $\nabla$ is the gradient with respect to the metric $\overbar g$. 
\end{itemize}

When $n>2$, we extend the metric $g$ on $I^n = [0, 1]^n$ to a complete metric $\overbar g$ on $\R^n$ by a similar induction argument. let us  briefly describe the construction below.

For each $k=0,1,\ldots,n$, let $\Lambda_k$ be the family of all subsets of $\{1,2,\cdots,n\}$ with $k$-elements. For $\lambda\in \Lambda_k$, set
$$X_\lambda=\big\{(x_1,x_2,\cdots,x_n)\in\R^n: \text{if $j\notin \lambda$, then $0\leq x_j\leq 1$} \big\}.$$
In particular, $X_\varnothing=[0,1]^n$ for $\lambda=\varnothing\in\Lambda_0$. For each $\lambda\in\Lambda_k$, the boundary of $X_\lambda$ has $(n-k)$ pairs of opposite faces, denoted by $\partial_{j\pm}X_\lambda$ for each $j\notin \lambda$. 

Now we shall extend the metric $g$ on  $X_\varnothing = [0,1]^n$ smoothly to a complete metric $\overbar g$ on $\R^n$ by induction on $k$  so that for each given $k\in\{0,1,\cdots,n-1\}$,  we have 
\begin{equation}\label{eq:dist}
\dist_{\overbar g}(\partial_{j-}X_\lambda,\partial_{j+}X_\lambda)\geq \ell_j-k\varepsilon/(n+1) \textup{ for } \forall\lambda\in\Lambda_k \textup{ and } \forall j\not\in\lambda,~
\end{equation} 
and the metric $\overbar g$ coincides with the Euclidean metric on $\mathbb R^n$  outside a compact set.
For $k=0$, this is trivial. Now assume that the metric has been constructed on $X_\lambda$ satisfying the condition \eqref{eq:dist}, for each $\lambda\in\Lambda_k$. For each $\mu\in\Lambda_{k+1}$, it suffices to extend the metric on the union  $ \cup_{\lambda\subset\mu}X_\lambda$ to $X_{\mu}$, where $\lambda$ runs through all $k$-element subsets of  $\mu$. Note that the condition  \eqref{eq:dist} is already satisfied on $\cup_{\lambda\subset\mu}X_\lambda$. Therefore such an extension of metric always exists on $X_{\mu}$ by the Whitney's extension theorem, a smooth partition of unity, and rescaling. By induction, we obtain a metric $\overbar g$ on $\R^n$ with the required properties. 

For each $i=1,2,\cdots, n$, set
$$E_i=\{(y_1,\cdots,y_n)\in\R^n:y_i=0 \}.$$
For each $x=(x_1,x_2,\cdots,x_n)\in\R^n$, we define
\begin{equation}\label{eq:phi}
\varphi_i(x)\coloneqq \sgn(x_i)\cdot \dist_{\overbar g}(x,E_i),
\end{equation}
where $\sgn$ is the sign function. 
 Set $\psi_i=\varphi_i-\ell_i/2$. For each $i=1,2, \cdots, n$, let $z_i$ be a smooth approximation of $\psi_i$ such that  $|z_i-\psi_i|\leq\varepsilon/(n+1)$ and
$\|\nabla z_i\|\leq 1+\varepsilon$. Here $\nabla$ is the gradient with respect to the metric $\overbar g$. 
To summarize, we have constructed a smooth metric $\overbar g$  and smooth functions $\{z_i\}_{i=1,2,\cdots,n}$ on $\R^n$ with the following properties.
\begin{itemize}
	\item The metric $\overbar g$ restricts to the metric $g$ on $[0,1]^n$, and coincides with the Euclidean metric on $\R^n$ outside a compact set. In particular, there exists a positive number $\sigma>0$ such that $\Sc_{\overbar g}\geq -\sigma$ on $\R^n$.
	\item For each $x\in\R^n$, if $|z_i(x)|\leq \ell_i/2-\varepsilon$ for $1\leq i\leq n$, then $x$ lies in the cube $[0,1]^n$.
	\item For each $C > 0$, the subset $\{x\in \R^n: z_1(x)^2+\cdots +z_n(x)^2\leq C\}$ is compact.
	\item $\|\nabla z_i\|\leq 1+\varepsilon$, where $\nabla$ is the gradient with respect to the metric $\overbar g$. 
\end{itemize}

With the constant $\varepsilon>0$ and $\sigma>0$ given above, there exist positive numbers $\delta>0$ and $\{\overbar\ell_i\}_{1\leq i\leq n}$ such that 
for all $1\leq i\leq n$, we have 
\begin{equation}
0<\overbar \ell_i  < \ell_i-2\varepsilon< \ell_i, 
\end{equation}
\begin{equation}\label{eq:1}
\frac{n^2}{4}-\sum_{i=1}^{n}\frac{\pi^2(1+\varepsilon)^2}{(\ell_i-2\varepsilon)^2}>\delta,
\end{equation}
and
\begin{equation}\label{eq:2}
\frac{-n\sigma}{4(n-1)}+
\frac{\pi^2(1+\varepsilon)^2}{(\ell_1-2\varepsilon)^2}\tan^2\frac{\pi\overbar\ell_1}{2(\ell_1-2\varepsilon)}
-\varepsilon-\sum_{i=2}^{n} \frac{\pi^2(1+\varepsilon)^2}{(\ell_i-2\varepsilon)^2}>\delta.
\end{equation}

Set 
$$r_i=\frac{\pi(1+\varepsilon)}{\ell_i-2\varepsilon}.$$
Fix a positive number $\ell$ such that $\ell>\ell_i/2+1$ for all $1\leq i\leq n$. We choose smooth functions $\{\xi_i\}_{1\leq i \leq n}$ on
 $\mathbb E^n$, the flat Euclidean space, such that
\begin{enumerate}[(1)]
	\item $0\leq (1+\varepsilon)\xi_i(x,y)\leq y^2+r_i^2$ for $\forall (x,y)\in\mathbb E$, 
	\item $0\leq \xi_i(x,y)\leq \varepsilon$ for $|x|\geq \ell_i/2$, where $\varepsilon$ is the constant given in line \eqref{eq:1}, 
	\item and 
	\[  \xi_i(x,y)= \begin{cases}
	\frac{y^2+r_i^2}{1+\varepsilon} & \textup{ if }  |x|\leq \overbar \ell_i/2,\\ 
	0 & \textup{ if }   \ell_i/2-\varepsilon\leq |x|\leq \ell, \\
	\varepsilon & \textup{ if } |x|\geq \ell +1.
	\end{cases}\]
\end{enumerate}

Let $f_i$ be the unique solution to the following differential equation with initial condition
\begin{equation}\label{eq:ode}
\begin{cases}
f'(x)=\xi_i(x,f(x)),\\
f(0)=0.
\end{cases}
\end{equation}
The solution always exists in a small neighborhood of $x=0$. 
By the comparison theorem, we have
\begin{equation}
f_i(x)\leq r_i\tan(\frac{r_i x}{1+\varepsilon}),
\end{equation}
since the function $h_i(x)= r_i\tan(\frac{r_i x}{1+\varepsilon})$ is the unique solution to the differential equation 
$$\begin{cases}
f'(x)=\frac{f(x)^2+r_i^2}{1+\varepsilon},\\
f(0)=0.
\end{cases}$$
Therefore the solution $f_i$ exists at least for $|x|<\ell_i/2$. On the other hand, by assumption, $\xi_i(x,y) = 0$ for all $\ell_i/2-\varepsilon\leq |x|\leq \ell$ and $0\leq \xi_i(x,y)\leq \varepsilon$ for $|x|\geq \ell_i/2$. It follows that the solution $f_i$ exists on the whole real line (cf. \cite[\S9, Theorem XIII]{MR1629775}). Roughly speaking, $f_i(x)$ is equal to the function  $h_i(x) = r_i\tan(\frac{r_ix}{1+\varepsilon})$ when $|x|$ is small, and is equal to the linear function $\varepsilon x \pm  c$ for some constant $c$ when  $|x|$ is large. 

Let $\cl_{n, 0}$ be the Clifford algebra generated by
$\{e_i\}_{1\leq i \leq n}$ subject to the relation $e_i^2=-1$ and $e_ie_j+e_je_i=0$ for $i\ne j$. Similarly, we define $\cl_{0, n}$ to be the Clifford algebra generated by $\{\hat e_i\}_{1\leq i \leq n}$ subject to the relation $\hat e_i^2=1$ and $\hat e_i\hat e_j+\hat e_j\hat e_i=0$ for $i\ne j$. Let $T\R^n\oplus \mathbb E^n$ be the direct sum of the tangent bundle of $(\R^n,\overbar g)$ and the trivial bundle $\mathbb E^n$, and
$$E\coloneqq S(T\R^n\oplus \mathbb E^n)$$ 
the spinor bundle of $T\R^n\oplus \mathbb E^n$. By construction, $E$ is equipped with the Clifford action of $\cl_{n,0}(T\R^n)\hotimes \cl_{0,n}(\mathbb E^n)$.

Let $D$ be the twisted Dirac operator on $L^2(\R^n,E)$, where  $L^2(\R^n,E)$ is the space  of $L^2$ sections of the bundle $E$.
In terms of local orthonormal frames $\{w_i\}_{1\leq i \leq n}$ of the tangent bundle $T\R^n$, the operator  $D$ can be locally expressed by 
\[ D = \sum_{i=1}^n w_i\nabla_{w_i}, \]
where $w_i$ denotes the Clifford action of $\cl_{n,0}(T\R^n)$ on $E$.

Consider the following Callias-type operator
$$B=D\hotimes 1+\sum_{i=1}^n 1\hotimes f_i(z_i)\hat e_i$$
where $\{\hat e_i\}$ is the flat orthonormal basis of $\mathbb E^n$. A direct computation shows that
\begin{equation}\label{eq:B2}
B^2=D^2\hotimes 1+\sum_{i=1}^n f_i(z_i)^2+\sum_{i=1}^n[D,f_i(z_i)]\hotimes \hat e_i.
\end{equation}
As $[D,f_i(z_i)]\hotimes \hat e_i$ is a bounded operator on $L^2(\R^n,E)$ for all $1\leq i\leq n$ and $\sum_{i=1}^n f_i(z_i)^2$ is a proper\footnote{A function $f\colon \R^n \to \R$ is said to be proper if the preimage $f^{-1}(K)$ of each compact subset $K$ is compact.} function on $\R^n$, it follows that $B$ is essentially self-adjoint and Fredholm \cite{MR1233861}. 

Clearly, there is  a smooth path of Riemannian metrics $g_t$ on $\R^n$ connecting the metric $\overbar g$ to the standard Euclidean metric $g_{o}$ on $\R^n$. We may assume that there exists a compact set $K$ such that each metric $g_t$ coincides with the standard Euclidean metric on $\R^n$ outside $K$. Furthermore, there is a homotopy of functions $\{\theta_{t, i}: 0\leq t\leq 1\}$ between $f_i$ and the standard  coordinate function $x_i$ of $(\mathbb R^n, g_{o})$ such that $\sum_{i=1}^n \theta_{t, i}^2$ is a proper function  and $\theta_{i,t}$ has uniformly bounded gradient for all $t\in [0, 1]$. 
It follows from the homotopy invariance of the Fredholm index that  $\ind(B)=\ind(B_o)$, where $\ind(B_o)$ is the index of the following Bott-Dirac operator
\begin{equation} \label{eq:bottdirac}
B_{o}=\sum_{i=1}^n\frac{\partial}{\partial x_i}e_i\hotimes 1+\sum_{i=1}^n 1\hotimes x_i\hat e_i,
\end{equation}
acting on $L^2(\mathbb R^n,S(\mathbb E^n\oplus\mathbb E^n))$.

We have 
$$B_{o}^2=-\sum_{i=1}^n\frac{\partial^2}{\partial x_i^2}+\sum_{i=1}^n x_i^2+\sum_{i=1}^n e_i\hotimes \hat e_i.$$
Observe that $\{e_i\hotimes\hat e_i\}_{1\leq i\leq n}$ is a commuting family of symmetries, i.e.,  $(e_i\hotimes\hat e_i)^2=1$ and 
$$(e_i\hotimes\hat e_i)(e_j\hotimes\hat e_j)=(e_j\hotimes\hat e_j)(e_i\hotimes\hat e_i)$$
for $i\ne j$.  Let $\Lambda$ be the set of all maps from  $\{1,2,\cdots,n\}$ to $ \{1,-1\}.$ The family of symmetries $\{e_i\hotimes\hat e_i\}_{1\leq i\leq n}$ decompose  
$S(\mathbb E^n\oplus\mathbb E^n)$
into $2^n$ subspaces $V_\lambda$ indexed by $\lambda\in \Lambda$, where
$$(e_i\hotimes\hat e_i)f=\lambda(i)f \textup{ for } \forall f\in V_\lambda.$$ 
As the rank of $S(\mathbb E^n\oplus\mathbb E^n)$ is equal to $2^n$, each $V_\lambda$ has rank $1$.

The operator $B_o^2$ also decomposes correspondingly into an orthogonal  direct sum of operators 
$$-\sum_{i=1}^n\frac{\partial^2}{\partial x_i^2}+\sum_{i=1}^n x_i^2+\sum_{i=1}^n \lambda(i)$$ 
each of which acts on $L^2(\R^n, V_\lambda)$.
It is known that the harmonic oscillator 
\[ -\frac{d^2}{dx^2}+x^2 \]
has a complete system of eigenvectors in $L^2(\R)$ with  the set of eigenvalues $\mathscr O=\{1,3,5,\ldots\}$.  Therefore the operator 
$$-\sum_{i=1}^n\frac{\partial^2}{\partial x_i^2}+\sum_{i=1}^n x_i^2$$
has a complete system of eigenvectors in $L^2(\R^n)$ with the set of eigenvalues 
$$\mathscr O_n\coloneqq \{k_1+k_2+\cdots+k_n: k_i\in\mathscr O \}.$$
It follows that  $B^2$ restricted on $L^2(\R^n,V_\lambda)$ is one-to-one and onto unless $\lambda$ is the constant map with $\lambda(i)=-1$ for all $i \in \{1,2,\cdots,n\}$. In particular, when $\lambda\equiv -1$, the kernel of $B^2$ is
$\{ \exp(-(x_1^2+\cdots +x_n^2)/2)\cdot v: v\in V_\lambda\},$
which is  of dimension $1$. As such $V_\lambda$ is located in the even part of $S(\mathbb E^n\oplus\mathbb E^n)$ with respect to its even-odd grading, the (graded) Fredholm index of $B$ is equal to $1$.

On the other hand, we shall show that $B$ is an invertible operator, hence has $\ind(B) = 0$,  which will lead to a contradiction. Indeed,  by the Lichnerowicz's formula, we have 
\[ D^2 = \nabla^\ast \nabla + \frac{\Sc_{\overbar g}}{4}.  \]
It follows from  the Cauchy–Schwarz inequality that 
\[  \langle D v, D v\rangle \leq n \langle \nabla v, \nabla v\rangle   \]
for all $v\in C_c^\infty(\mathbb R^n, E)$. Therefore, we have 
\begin{equation}
D^2\geq \frac{\Sc_{\overbar g}}{4}\cdot \frac{n}{n-1} \textup{ on } C_c^\infty(\mathbb R^n, E).
\end{equation}
By line \eqref{eq:B2}, we have 
\begin{align*}
B^2&\geq D^2+\sum_{i=1}^nf_i(z_i)^2-\sum_{i=1}^n\|[D,f_i(z_i)]\|\\
&\geq \frac{\Sc_{\overbar g}}{4}\cdot \frac{n}{n-1}+\sum_{i=1}^nf_i(z_i)^2-\sum_{i=1}^nf'_i(z_i)\|\nabla z_i\|\\
&\geq \frac{\Sc_{\overbar g}}{4}\cdot \frac{n}{n-1}+\sum_{i=1}^nf_i(z_i)^2-\sum_{i=1}^nf'_i(z_i)(1+\varepsilon).
\end{align*}

\textbf{Case (1): inside the cube.}  By assumption, we have  $\Sc_{\overbar g}\geq n(n-1)$ inside the cube $I^n = [0, 1]^n$. Hence
\begin{align*}
B^2&\geq \frac{n^2}{4}+\sum_{i=1}^n f_i(z_i)^2-\sum_{i=1}^nf'_i(z_i)(1+\varepsilon)\\
&=\frac{n^2}{4}+\sum_{i=1}^n f_i(z_i)^2-\sum_{i=1}^n\xi_i(z_i,f_i(z_i))(1+\varepsilon)\\
&\geq \frac{n^2}{4}-\sum_{i=1}^nr_i^2=
\frac{n^2}{4}-\sum_{i=1}^n\frac{\pi^2(1+\varepsilon)^2}{(\ell_i-2\varepsilon)^2}>
\delta
\end{align*}
where the positivity is deduced from line \eqref{eq:1}.

\textbf{Case (2): outside the cube.} In this case, we see that at least one of the $|z_i|$'s is $\geq \ell_i/2-\varepsilon$. Without loss of generality, we assume that $|z_1|\geq\ell_1/2-\varepsilon>\overbar \ell_1/2$. Note that $f_1(x)=r_1\tan(\frac{r_1x}{1+\varepsilon})$ for $|x|\leq \overbar \ell_1/2$. Since $f_1$ is an increasing function, we have 
$$|f_1(x)|\geq r_1\tan\frac{r_1\overbar \ell_1}{2(1+\varepsilon)}=
\frac{\pi(1+\varepsilon)}{\ell_1-2\varepsilon}\tan\frac{\pi\overbar\ell_1}{2(\ell_1-2\varepsilon)}
\textup{ for all } |x|\geq \ell_1/2-\varepsilon.$$
Furthermore, since $|z_1|\geq \ell_1/2-\varepsilon$, we have $|f_1'(z_1)|\leq \varepsilon$. 
Therefore, we see that
\begin{align*}
B^2&\geq \frac{-n\sigma}{4(n-1)}+f_1(z_1)^2-\varepsilon-\sum_{i=2}^{n} r_i^2\\
&\geq \frac{-n\sigma}{4(n-1)}+
\frac{\pi^2(1+\varepsilon)^2}{(\ell_1-2\varepsilon)^2}\tan^2\frac{\pi\overbar\ell_1}{2(\ell_1-2\varepsilon)}
-\varepsilon-\sum_{i=2}^{n} r_i^2>\delta,
\end{align*}
where the positivity is deduced from line \eqref{eq:2}. 

To summarize, we have proved that $B^2 >\delta$ for some positive constant $\delta>0$. It follows that $B$ is invertible, hence $\ind(B) = 0$, which contradicts the fact that $\ind(B)=1$. This finishes Part 1 of the proof.

\proofpart{2}{Strict inequality}
	
	Given a Riemannian metric on the cube $[0,1]^n$ with scalar curvature $\geq n(n-1)$, we assume to the contrary that
		  \begin{equation}\label{eq:cube>=}
	\sum_{i=1}^{n}\frac{1}{\ell_i^2} =  \frac{n^2}{4\pi^2}. 
	\end{equation}
	
	Let  $z_i$ and $f_i$ be functions similar to those defined in Part 1, for $i=1,2,\cdots,n$. We shall specify more precisely what properties $z_i$ and $f_i$ need to satisfy,  after we set up a relevant index problem and its associated Bochner-Lichnerowicz type formula. For the moment, we only assume each $z_i$ is a real-valued smooth function on $I^n$ and $f_i$ is a real-valued smooth function of $z_i$.  
		
	Consider the bundle 
	$$E\coloneqq S(TI^n\oplus \mathbb E^n)$$
	 over the cube $I^n$, similarly to Part 1 . Set
\begin{equation}\label{eq:operatorB}
		B=D\hotimes 1+\sum_{i=1}^n 1\hotimes f_i(z_i)\hat e_i.
\end{equation}
	We shall introduce a boundary condition $\mathfrak  B$ at the codimension one faces of $I^n$ for  sections of $E$ over $\partial I^n$. By assumption, the dihedral angles of $I^n$ (under the metric $g$) are $ < \pi$, since the corner structure is non-degenerate.
	
Since the dihedral angles of $I^n$ (under the metric $g$) are $<\pi$,  it is not difficult to see that there exists a manifold with corners $M$ in $\R^n$ and a smooth map $\Phi\colon I^n\to M$ such that
	\begin{enumerate}
		\item $\Phi$ preserves the corner structures,
		\item $0< \theta(x) <  \theta_M(\Phi(x))< \pi$ for all $x$ in each  codimension two face of $I^n$, where $\theta(x)$ is the dihedral angle of $(I^n, g)$ at the  point $x$   and $\theta_M(\Phi(x))$ is the corresponding dihedral angle of $M$ at $\Phi(x)$. 
		\item For each $x$ in a codimension one face $F$ of $M$, we have 
		\[ \langle \hat\nu_x,\hat e_i\rangle < 0   \textup{ if } F =  \Phi(\partial_{i,+}), \textup{ and } \langle \hat\nu_x,\hat e_i\rangle > 0 \textup{ if } F = \Phi(\partial_{i,-}), \] where $\hat\nu_x$ is the unit inner normal vector field of $F$ at $x$ and $\hat e_i$ is the $i$-th standard basis vector\footnote{Note that we have used $\hat e_i$ to denote both $i$-th standard basis vector of $\mathbb R^n$ and the associated Clifford multiplication on $\cl_{0, n}$.  It should be clear from the context which one we refer to. } of $\mathbb E^n$.
	\end{enumerate}

	For example, we may slightly perturb the unit sphere in $\mathbb R^n$ and choose $M$ to be the compact region enclosed by this slightly perturbed sphere.  For convenience,\footnote{In the present paper, the convexity of $\partial M$ is only assumed for convenience, and is not really essential for the proof.} we may further  assume that each codimension one face of $M$ is convex (that is, has   non-negative second fundamental form).  Now the boundary condition $\mathfrak B$ on each codimension one face $\overbar F$ of $I^n$  is given by
	\begin{equation}\label{eq:boundaryCondition}
		\big(\nu_x \hotimes \hat\nu_{\Phi(x)} \big)\varphi(x)=-\varphi(x)
	\end{equation}
	for all $x\in \overbar F$ and all sections $\varphi$ of $E$, where $\nu$ is the Clifford multiplication given by the unit inner normal vector field $\nu_x$ of $\overbar F$ at $x$, and $\hat\nu_{\Phi(x)}$ is the Clifford multiplication given by the unit inner normal vector field $\hat\nu_{\Phi(x)}$ of $\Phi(\overbar F)$ at $\Phi(x)$.
	
	Note that the operator $D$ is a twisted Dirac operator to which the index theoretic techniques   developed in  \cite{Wang:2021tq} for manifolds with polyhedral boundary  apply. As $B$ only differs from $D$ by a bounded term,  the same techniques from \cite{Wang:2021tq} imply the following properties of $B$. 
	\begin{enumerate}
		\item The operator $B$ with the boundary condition $\mathfrak B$ is an essentially self-adjoint Fredholm operator with domain $H^1(I^n, E; \mathfrak B)$   consisting of Sobolev $H^1$-sections satisfying the boundary condition $\mathfrak  B$. 
		\item Furthermore, the Fredholm index of $B$ is equal to $1$. In particular, there exists a nonzero element $\varphi$ in $H^1(I^n, E; \mathfrak B)$  such that $B\varphi = 0$. 
	\end{enumerate}

	Let $\psi$ be any smooth section of $E$ over $I^n$ that vanishes near all codimension two faces. We have
	\begin{align}
		\int_{I^n}\langle B\psi,B\psi\rangle=&\int_{I^n} |D\psi|^2+
		\int_{I^n}\Big \langle \big(\sum_{i=1}^n f_i(z_i)^2+\sum_{i=1}^n[D,f_i(z_i)]\hotimes \hat e_i\big)\psi,\, \psi\Big \rangle \notag\\
		& + \int_{\partial I^n}\sum_{i=1}^n \langle f_i(z_i)(\nu\hotimes\hat e_i)\psi,\psi\rangle \label{eq:bochner}
	\end{align}
	Let $\mathcal P$ be the Penrose operator defined as follows: 
	\begin{equation}\label{eq:Penrose}
		\mathcal P_V\psi\coloneqq \nabla_V\psi+\frac 1 n V\cdot D\psi
	\end{equation}
	for $V\in TI^n$. We have the following identity (\cite[Section 5.2]{MR3410545}):
	 \begin{equation}\label{eq:P2+D2}
	 	|\nabla\psi|^2=|\mathcal P\psi|^2+\frac 1 n |D\psi|^2.
	 \end{equation}
	 Note that
	 \begin{align*}
	 	\int_{I^n} |D\psi |^2=&\int_{I^n}\langle D^2\psi,\psi\rangle+\int_{\partial I^n} \langle \nu\cdot D\psi,\psi\rangle\\
	 	=&\int_{I^n}\frac{\Sc_g}{4}|\psi|^2+\int \langle\nabla^*\nabla\psi,\psi\rangle+\int_{\partial I^n} \langle \nu\cdot D\psi,\psi\rangle\\
	 	=&\int_{I^n}\big(\frac{\Sc_g}{4}|\psi|^2+ |\nabla\psi|^2\big)+\int_{\partial I^n} \langle (\nu\cdot D+\nabla_\nu)\psi,\psi\rangle,
	 \end{align*}
 which together with Equation \eqref{eq:P2+D2} implies that 
 $$\int_{I^n} |D\psi|^2=\frac{n}{n-1}\int_{I^n}(\frac{\Sc_g}{4}|\psi|^2+|\mathcal P\psi|^2)+\frac{n}{n-1}\int_{\partial I^n} \langle (\nu\cdot D+\nabla_\nu)\psi,\psi\rangle.$$
On each codimension one face of $I^n$, we have that
 $$\nu\cdot D+\nabla_\nu=\sum_{j=1}^{n-1}e_j\nabla_{e_j},$$
 where $\{e_j\}$ is a local orthonormal frame of the tangent bundle of the face. In particular, 
 $$\mathscr H = [\nu\cdot D+\nabla_\nu,\,   \nu\hotimes\hat \nu ]\coloneqq (\nu\cdot D+\nabla_\nu)(\nu\hotimes\hat \nu)+(\nu\hotimes\hat \nu)(\nu\cdot D+\nabla_\nu)$$
 is a bounded endomorphism acting on the bundle $E$ over this face.\footnote{In fact, under the assumption that the codimension one faces of $M$ are convex,  for each $x$ in a codimension one face $\overbar F$ of $I^n$, a local computation (cf. \cite[Section 2]{Wang:2021tq}) shows that 
 $$\mathscr H(x)\geq  \frac 1 2 H(x)-\frac 1 2 H_{\partial M}(\Phi(x))\cdot \|d\Phi\|_x,$$
 where $H$ (resp. $H_{\partial M}$) is the mean curvatures of $\overbar F$ (resp. $\Phi(\overbar F)$) and $\|d\Phi\|_x$ is the norm of the map $d\Phi\colon T_xI^n \to T_{\Phi(x)}M$. Such an explicit lower bound of $\mathscr H$ is important in some geometric applications. However, for the present paper,  we actually only need to know $\mathscr H$ is (uniformly) bounded from below on the codimension one faces of $I^n$.  }
 
 Now we consider  the term $\langle (\nu\hotimes\hat e_i)\psi,\psi\rangle$ from line \eqref{eq:bochner}. Recall that $\{\hat e_i\}$ is the standard orthonormal basis of $\mathbb R^n$. If $\hat \nu$ is  the unit inner normal vector $ \hat \nu_y$ of a codimension one face of $M$ at $y$, then  we have the decomposition 
 \[ \hat e_i=a_i\hat\nu+ \hat w_i, \] where $a_i = \langle\hat \nu, \hat e_i\rangle$  and  $\hat w_i$ is a vector that is orthogonal to $\hat\nu$. It follows that 
 $$\langle (\nu\hotimes\hat e_i)\psi,\psi\rangle=a_i\langle (\nu\hotimes\hat\nu) \psi,\psi\rangle+\langle (\nu\hotimes\hat w_i)\psi,\psi\rangle.$$ 
 Now suppose $\psi$ satisfies the boundary condition $\mathfrak B$ from line \eqref{eq:boundaryCondition}. Then we have 
 \[  a_i\langle (\nu\hotimes\hat\nu) \psi,\psi\rangle  = a_i|\psi|^2 \ \textup{ and } \  \langle (\nu\hotimes\hat w_i)\psi,\psi\rangle = 0, \]
 where the latter follows from that fact that  $\nu\hotimes\hat w_i$ anti-commutes with $\nu\hotimes\hat\nu$. To summarize,
we have that 
		\begin{equation}\label{eq:Greenformula}
		\begin{split}
			&\int_{I^n}\langle B\psi,B\psi\rangle\\
			=&\int_{I^n} |\mathcal P\mathcal \psi|^2+
		\int_{I^n}\Big\langle \Big(\frac{n}{n-1}\frac{\Sc_g}{4}+\sum_{i=1}^n f_i(z_i)^2+\sum_{i=1}^n[D,f_i(z_i)]\hotimes \hat e_i\Big)\psi,\psi\Big\rangle\\
		&+\frac{n}{n-1}\int_{\partial I^n}\langle\mathscr H\psi,\psi\rangle+ \int_{\partial I^n} \big(\sum_{i=1}^n \langle\hat\nu, \hat e_i\rangle f_i(z_i)\big)|\psi|^2
		\end{split}
	\end{equation}
	for all smooth sections $\psi\in C_0^\infty(I^n, E; \mathfrak B)$, where $C_0^\infty(I^n, E; \mathfrak B)$ is the space of all smooth sections of $E$ over $I^n$ that vanish near all codimension two faces of $I^n$ and satisfy the boundary condition $\mathfrak B$ at all codimension one faces of $I^n$.  Note that the completion of $C_0^\infty(I^n, E; \mathfrak B)$ with respect to the Sobolev $H^1$ norm is $H^1(I^n, E; \mathfrak B)$, since removing a subspace of codimension at least two does not affect Sobolev $H^1$ space. As  both sides of Equation \eqref{eq:Greenformula} are continuous with respect to the $H^1$-norm, it follows that  Equation  \eqref{eq:Greenformula} also holds for all $\psi\in H^1(I^n, E; \mathfrak B)$. 

Now we shall describe the properties that  $f_i$ and $z_i$ need to satisfy. Set $\ell_i=\dist_g(\partial_{i-},\partial_{i+})$. We first define
\begin{equation}
	f_i(x)=-\frac{\pi}{\ell_i}\tan(\frac{\pi x}{\ell_i}),
\end{equation}
for $x\in(-\ell_i/2,\ell_i/2)$, which satisfies that
\begin{equation}\label{eq:odeForf}
	\frac{\pi^2}{\ell_i^2}+f_i^2-|f_i'|=0.
\end{equation}
Note that $f_i(x)\to\mp\infty$ as $x\to\pm\ell_i/2$. 
For each $1\leq i \leq n$, consider the function
\[ \max\{\ell_i,\dist(\cdot,\partial_{i-})\}-\ell_i/2\]
on $I^n$. By modifying it slightly on  $I-\partial_{i\pm}$, we see that there exists a function $\widetilde z_i$ on $I^n$ such that
\begin{itemize}
	\item  $\widetilde z_i$ is $1$-Lipschitz on $(I^n,g)$,
	\item $\widetilde z_i=\pm \ell_i/2$ on $\partial_{i\pm}$, and
	\item $|\widetilde z_i(x)|<\ell_i/2$ if $x$ is not in $\partial_{i\pm}$.
\end{itemize}
At the moment, $\widetilde z_i$ is not necessarily smooth on $I^n$.   But in any case, the function  $f_i(\widetilde z_i)$ is well-defined everywhere in the interior of $I^n$. 

By our choice of $M$,  we have that $a_i \coloneqq  \langle \hat v, \hat e_i\rangle $ is strictly positive (resp. negative) on $\Phi(\partial_{i-})$ (resp. $\Phi(\partial_{i+})$). Therefore we have 
\begin{equation}\label{eq:ithLimitBoundary}
	a_i f_i(\widetilde z_i(x)) \textup{ uniformly goes  to } +\infty, \textup{ as $x$ approaches } \partial_{i\pm}.
\end{equation}
We have the following claim. 
\begin{claim}\label{claim:limitboundary} 
		$\sum_{i=1}^n a_i f_i(\widetilde z_i(x))$ uniformly goes  to $+\infty$, as $x$ approaches $\partial I^n.$
\end{claim}	
Suppose $y$ is a point in the interior of the intersection of $k$ codimension one faces. For example,  assume  $y$ lies in the interior of  $\partial_{1+}\cap \cdots \cap \partial_{k+}$. By line \eqref{eq:ithLimitBoundary}, there exists a small neighborhood $U_y$ of $y$ in $N$ such that the functions  $a_1 f_1(\widetilde z_1),\cdots,a_k f_k(\widetilde z_k)$ are all positive on $U_y$, while the functions 
\[ a_{k+1} f_{k+1}(\widetilde z_{k+1}),\cdots,a_n f_n(\widetilde z_n)\] are uniformly bounded on $U_y$. Therefore $\sum_{i=1}^n a_i f_i(\widetilde z_i(x)) $  goes to $+\infty$ uniformly when $x$ approaches to $\partial I^n$, as long as $x\in U_y$. Now Claim \ref{claim:limitboundary} follows by the compactness of $\partial I^n$.

By applying Claim \ref{claim:limitboundary}, it is not difficult to see that we can slightly perturb $\widetilde z_i$ to obtain functions $z_i$ with the following properties.
\begin{enumerate}
	\item Each $z_i$ is smooth and $1$-Lipschitz on $(I^n,g)$,
	\item $|z_i(x)|<\ell_i/2$ for all $x\in I^n$, and
	\item for all $x\in \partial_{j\pm}$ and all $1\leq j\leq n$,  we have 
	\begin{equation}\label{eq:boundary>}
		\sum_{i=1}^n \langle \hat v, \hat e_i\rangle  f_i(z_i(x))>\frac{n}{n-1}\|\mathscr H\|
	\end{equation}
\end{enumerate}
Note that
\begin{equation}\label{eq:interior>=0}
	\begin{split}
		&\frac{n}{n-1}\frac{\Sc_g}{4}+\sum_{i=1}^n f_i(z_i)^2+\sum_{i=1}^n[D,f_i(z_i)]\hotimes \hat e_i\\
		\geq&\frac{n^2}{4}+\sum_{i=1}^n f_i(z_i)^2-\sum_{i=1}^n|f_i'(z_i)|\\
		=&\frac{n^2}{4}-\sum_{i=1}^n\frac{\pi^2}{\ell_i^2}=0
	\end{split}
\end{equation}
where the first inequality is because the scalar curvature is $\geq n(n-1)$ and each $z_i$ is $1$-Lipschitz,  the second equality follows from line \eqref{eq:odeForf}, and the last equality is due to   the assumption \eqref{eq:cube>=}.

Recall that the operator $B$ from  line \eqref{eq:operatorB} (subject to the boundary condition $\mathfrak B$) is an essentially self-adjoint Fredholm operator with domain $H^1(I^n, E; \mathfrak B)$. Moreover, the Fredholm index of $B$ is  $1$. So there exists a  non-zero element $\varphi \in H^1(I^n, E; \mathfrak B)$ such that $B\varphi=0$. We plug $\varphi$ into  Equation \eqref{eq:Greenformula}, it follows from line \eqref{eq:boundary>} and \eqref{eq:interior>=0} that
\begin{itemize}
	\item $\mathcal P\varphi=0$ almost everywhere on $I^n$, and
	\item $\varphi=0$ almost everywhere on $\partial I^n$.
\end{itemize}
Now $\mathcal P\varphi=0$  and $B\varphi=0$ imply that 
\begin{equation}\label{eq:odeForphi}
	\nabla_V\varphi=\frac 1 n  \left(\sum_{i=1}^n V\hotimes f_i(z_i)\hat e_i\right)\varphi,
\end{equation}
for all $V\in TI^n$, 
which is a first order linear ordinary differential equation of $\varphi$. It follows that $\varphi$ is smooth on $I^n$. On the other hand, we already know $\varphi=0$ almost everywhere on $\partial I^n$. In particular, $\varphi(y_0) = 0 $ for some point  $y_0$ in $\partial I^n$. Now for any $x\in I^n$, we choose a smooth path $\gamma_x$  connecting $x$ to $y_0$.  By the uniqueness of solutions to  the linear differential equation \eqref{eq:odeForphi} along the path $\gamma_x$, we see that $\varphi\equiv 0$ everywhere on $I^n$. But this contradicts the fact that $\varphi$ is a non-zero element in $H^1(I^n, E; \mathfrak B)$. This finishes the proof of Part 2. 
\end{proof}

\begin{remark}\label{rm:anglePi}
In fact, by combining  the approximation argument from \cite[Lemma 5.3]{Wang:2021tq} and \cite[Lemma 2.10]{WangXie2022}, the above proof for Part 2 of  Theorem \ref{thm:cube} also applies to the case where the dihedral angles of $g$ are allowed to be $\leq \pi$ instead of $< \pi$.  Here is a sketch. First, let $M$ be the unit round ball in $\mathbb R^n$. We construct similar functions $f_i(z_i)$ as in the above proof for Part 2 of  Theorem \ref{thm:cube}. However, we do not compute the index of $B$ with boundary condition $\mathfrak B$. Instead, we  approximate the metric $g$ on $I^n$ by a sequence of metrics $g_k$  such that the dihedral angles of $g_k$ are strictly less than $\pi$. The corresponding twisted Dirac operator $B_k$ on $(I^n, g_k)$ subject to the corresponding boundary condition $\mathfrak B_k$ is Fredholm with index equal to $1$, hence admits nontrivial element $\varphi_k\in  H^1(I^n_{g_k}, E; \mathfrak B_k)$ such that $B_k\varphi_k=0$.  By the non-negativity of each term in line \eqref{eq:Greenformula}, it follows from  \cite[Lemma 5.3]{Wang:2021tq} and \cite[Lemma 2.10]{WangXie2022} that the sequence $\{\varphi_k\}$ has a subsequence convergent to a nonzero element $\varphi$ in $H^1(I^n_{g}, E; \mathfrak B)$ such that $B\varphi =0$. Now the rest of the proof proceeds exactly the same way as the proof for Part 2 of Theorem \ref{thm:cube}.
\end{remark}

\begin{remark}\label{rk:optimal}
	The constant $4\pi^2$ in line \eqref{eq:cube} of Theorem \ref{thm:cube} is optimal. In fact, for any $\varepsilon>0$, there is a Riemannian metric $g$ on $[0,1]^n$ with scalar curvature equal to $n(n-1)$, and
	\begin{equation}\label{eq:optimal}
	\sum_{i=1}^n\frac{1}{\ell_i^2}\leq \frac{n^2}{4\pi^2}+\varepsilon.
	\end{equation}
	For any $l<\pi/n$, consider the warped metric 
	$$g=dx_1^2+\varphi(x_1)^2\cdot \sum_{i=2}^n R^2dx_i^2$$
	on $[-l,l]\times [0,1]^{n-1}$, where $R$ is a positive constant and
	$$\varphi(x_1)=	\big(\cos\frac{nx_1}{2}\big)^{2/n}.$$
	As computed in \cite[page 653]{MR3816521}, the scalar curvature of $g$ is given by
	$$Sc_{g}=-2(n-1)\frac{\varphi''}{\varphi}-(n-1)(n-2)\frac{\varphi'^2}{\varphi^2}=n(n-1).$$
	Note that there is $C_l>0$ (depending on $l$) such that
	$$	\sum_{i=1}^{n}\frac{1}{\ell_i^2}\leq  \frac{1}{(2l)^2}+\frac{C_l}{R^2}.$$
	The optimality \eqref{eq:optimal} follows by first choosing $l$ to be close to $\pi/n$ and then choosing $R$ to be large enough.
\end{remark}

\section{Proof of Gromov's $\square^{n-m}$ inequality}
In this section, we prove Theorem \ref{thm:cubespin}. We assume that the constant  $k$ in Theorem 1.2 to be $n(n-1)$, while the general case follows by rescaling the metric.
\begin{proof}[Proof of Theorem \ref{thm:cubespin}]
	We first prove part 1 of Theorem \ref{thm:cubespin}. 
	\proofpart{1}{$\geq$ inequality}
	Let us assume to the contrary that 	
	\[ \sum_{i=1}^m \frac{1}{\ell_i^2} < \frac{n^2}{4\pi^2}. \]

	We first show that the general case where $\iota\colon \pi_1(Y_\pitchfork) \to \pi_1(X)$ is injective  can be reduced to   the case where $\iota\colon \pi_1(Y_\pitchfork) \to \pi_1(X)$ is split injective.\footnote{We say $\iota\colon \pi_1(Y_\pitchfork) \to \pi_1(X)$ is split injective if there exists a group homomorphism $ \varpi\colon \pi_1(X) \to \pi_1(Y_\pitchfork) $ such that $\varpi\circ \iota = \mathbf{1}$, where $\mathbf{1}$ is the identity morphism of $\pi_1(Y_\pitchfork)$.} This reduction step was shown in \cite{Xie:2021tm}. We repeat the proof here for the convenience of the reader. Let $X_{u}$ be the universal cover of $X$. Since by assumption $\iota\colon \pi_1(Y_\pitchfork) \to \pi_1(X)$ is injective, we can view $\Gamma = \pi_1(Y_\pitchfork)$ as a subgroup of $\pi_1(X)$.  Let $X_\Gamma = X_u/ \Gamma$ be the covering space of $X$ corresponding to the subgroup $\Gamma\subset \pi_1(X)$. Then the inverse image of $Y_\pitchfork$ under the projection $p\colon X_\Gamma \to X$ is a disjoint union of covering spaces of $Y_\pitchfork$, at least one of which is a diffeomorphic copy of $Y_\pitchfork$. Fix such a copy of $Y_\pitchfork$ in $X_\Gamma$ and denote it by $\widehat Y_\pitchfork$.  Roughly speaking, the space $X_\Gamma$ equipped with the lifted Riemannian metric from $X$ could serve as a replacement of the original space $X$, except that $X_\Gamma$ is not compact in general. To remedy this, we shall choose a ``fundamental domain" around $\widehat Y_\pitchfork$ in $X_\Gamma$ as follows. 
	
	By assumption, $Y_\pitchfork \subset X$ is the transversal intersection of $m$ orientable hypersurfaces $Y_i \subset X$. Let $r_i$ be the distance function\footnote{To be precise, let $r_i$ be a smooth approximation of the distance function from $\partial_{i-}$.} from $\partial_{i-}$, that is $r_i(x) = \dist(x, \partial_{i-})$. Without loss of generality,  we can assume $Y_i = r_i^{-1}(a_i)$ for some regular value $a_i\in [0, \ell_i]$.  Let $Y^\Gamma_i = p^{-1}(Y_i)$ be the inverse image of $Y_i$ in $X_\Gamma$. Denote by $\overbar{r}_i$ the lift of $r_i$ from $X$ to $X_\Gamma$. Let $\nabla \overbar{r}_i$ be the gradient vector field associated to $\overbar{r}_i$. A point $x\in X_\Gamma$ said to be \emph{permissible} if there exist a number $s\geq 0$ and a piecewise smooth curve $c\colon [0, s] \to X_\Gamma$ satisfying the following conditions: 
	\begin{enumerate}[(i)]
		\item $c(0) \in \widehat Y_{\pitchfork}$ and $c(s) = x$; 
		\item there is a subdivision of $[0, s]$ into finitely many subintervals $\{[t_{k}, t_{k+1}]\}$ such that,  on each subinterval $[t_{k}, t_{k+1}]$, the curve $c$ is either an integral curve or a reversed integral curve\footnote{By definition, an integral curve of a vector field is a curve whose tangent vector coincides with the given vector field at every point of the curve. A reversed integral curve is an integral curve with the reversed parametrization, that is,  the tangent vector field of a reserved integral curve coincides  with the negative of the given vector field at every point of the curve. } of the gradient vector field $\nabla \overbar{r}_{j_k}$ for some $1\leq j_k\leq m$, where we require $j_{k}$'s to be all distinct from each other;
		\item furthermore, when $c$ is an integral curve of the gradient vector field $\nabla \overbar{r}_{j_k}$ on the subinterval $[t_{k}, t_{k+1}]$, we require the length of $c|_{[t_{k}, t_{k+1}]}$ to be less than or equal to $(\ell_{j_k} - a_{j_k} -\frac{\varepsilon}{4})$; and when $c$ is a reversed integral curve of the gradient vector field $\nabla \overbar{r}_{j_k}$ on the subinterval $[t_{k}, t_{k+1}]$, we require the length of $c|_{[t_{k}, t_{k+1}]}$ to be less than or equal to $(a_{j_k} -\frac{\varepsilon}{4})$.
	\end{enumerate} 
	
	Let $T$ be the set of all permissible  points. Now $T$ may not be a manifold with corners. To fix this, we choose an open cover $\mathscr U = \{U_\alpha\}_{\alpha\in \Lambda}$ of $T$ by geodesically convex metric balls of sufficiently small radius $\delta > 0 $. Now take the union of members of $\mathscr U = \{U_\alpha\}_{\alpha\in \Lambda}$ that do not intersect the boundary $\partial T$ of $T$, and  denote by $Z$ the closure of the resulting subset. Then $Z$ is a manifold with corners which,   together with the subspace $\widehat Y_\pitchfork\subset Z$,  satisfies all the conditions of the theorem, provided that $	\varepsilon$ and $\delta$ are chosen to be sufficiently small. In particular, the intersection $Y_i^\Gamma\cap Z$ of each hypersurface $Y_i^\Gamma$ with $Z$ gives a hypersurface of $Z$.  The transerval intersection of the resulting hypersurfaces is precisely $\widehat Y_\pitchfork \subset Z$.   Furthermore, note that the isomorphism $\Gamma = \pi_1(Y^\Gamma_\pitchfork) \to \pi_1(X^\Gamma) = \Gamma$ factors as   the composition $\pi_1(Y^\Gamma_\pitchfork) \to \pi_1(Z) \to \pi_1(X^\Gamma)$, where the morphisms $\pi_1(Y^\Gamma_\pitchfork) \to \pi_1(Z)$ and  $\pi_1(Z) \to \pi_1(X^\Gamma)$ are induced by the obvious inclusions of spaces. It follows that   $\pi_1(Y^\Gamma_\pitchfork) \to \pi_1(Z)$ is a split injection. 
	Therefore, without loss of generality, it suffices to prove the theorem under the additional assumption that $\iota \colon \pi_1(Y_\pitchfork) \to \pi_1(X)$ is a split injection.

	From now on, let us assume $\iota \colon \pi_1(Y_\pitchfork) \to \pi_1(X)$ is a split injection with a
	splitting morphism $\overbar\omega: \pi_1 (X) \to \pi_1(Y_\pitchfork) =\Gamma$. Let $\widetilde X$
	be the Galois $\Gamma$-covering space
	determined by  $\overbar\omega: \pi_1 (X) \to \Gamma$. In particular, the restriction of the covering map $\widetilde X\to X$ on $Y_\pitchfork$
	gives the universal covering space of $Y_\pitchfork$. 

	For any sufficiently small $\varepsilon>0$ and for each $1\leq i\leq m$, there exists a real-valued smooth function $\varphi_i$ on $X$
	such that (cf. \cite[Proposition 2.1]{MR532376})
	\begin{enumerate}
		\item $\|\nabla \varphi_i \|\leq 1+\varepsilon$,
		\item and $\varphi_i=0$ on $\partial_{i-}$ and $\varphi_i\geq \ell_i$ on $\partial_{i+}$.
	\end{enumerate}
	Set $\psi_i=\varphi_i-\ell_i/2$.
	We may assume that $Y_\pitchfork=\cap_{i=1}^m \psi_i^{-1}(0)$. Since $Y_\pitchfork$ has a trivial normal bundle, there is $c>0$ such that $\cap_{i=1}^m \psi_i^{-1}([-c,c])$ is diffeomorphic to $Y_\pitchfork\times I^m = Y_\pitchfork\times [0,1]^m$. This shows that $\partial X$ is cobordant to $Y_\pitchfork\times \partial I^m$ via a manifold $W$. Set
	$$Z=X\cup_{\partial X}W^{\text{op}}\cup_{Y_\pitchfork\times \partial I^m} \big(Y_\pitchfork\times \partial I^m \times \R_{\geq 0}\big) ,$$
	where $W^{\text{op}}$ is the manifold $W$ with reversed orientation and $ \R_{\geq 0} = [0, \infty)$.
	
	Fix a metric $g_{Y_\pitchfork}$ on $Y_\pitchfork$. We shall construct a metric $\overbar g$ on $Z$ that extends $g$ on $X$ and  is equal to the product metric $g_{Y_\pitchfork}+dx_1^2+\cdots +dx_{m}^2$ on $Y_\pitchfork \times \mathbb R^m$ outside a compact subset of $Z$.
	
	 \begin{figure}[h]
		\centering
	\scalebox{0.9}{	\begin{tikzpicture}
		\draw [dashed](0.3,0.3) -- (0.3,-0.3) -- (-0.3,-0.3) -- (-0.3,0.3) -- (0.3,0.3);
	  	\draw (0,0) circle ({2*cos(45)});
	  	\filldraw (0,-1.7) node {$X$};
		\draw [dashed] (0.3,0.3) .. controls (0.5,0.7).. (1,1);
		\draw [dashed] (0.3,-0.3) .. controls (0.5,-0.7).. (1,-1);
		\draw [dashed] (-0.3,0.3) .. controls (-0.5,0.7).. (-1,1);
		\draw [dashed] (-0.3,-0.3) .. controls (-0.5,-0.7).. (-1,-1);
		\draw (1,1) .. controls (1.5,2) .. (4,2) -- (6,2);
		\draw (1,-1) .. controls (1.5,-2) .. (4,-2) -- (6,-2);
		\draw (-1,1) .. controls (-1.5,2) .. (-4,2) -- (-6,2);
		\draw (-1,-1) .. controls (-1.5,-2) .. (-4,-2) -- (-6,-2);
		\draw [dashed](4,2) -- (4,-2);
		\draw [dashed](-4,-2) -- (-4,2);
		\filldraw (0.9,0) node {$W_{1+}$};
		\filldraw (-0.9,0) node {$W_{1-}$};
		\filldraw (2.5,0) node {$W_{1+}^\text{op}$};
		\filldraw (-2.5,0) node {$W_{1-}^\text{op}$};
		\filldraw (6,0) node {$Y_\pitchfork\times[0,+\infty) \times[0,1]$};
		\filldraw (-6,0) node {$Y_\pitchfork\times(-\infty,0] \times[0,1] $};
		\filldraw (0,1) node {$W_{2+}$};
		\filldraw (0,-1) node {$W_{2-}$};
		\end{tikzpicture}
	}
		\caption{Extend the metric from $X$ to $W_2$.}
		\label{fig:W2}
	\end{figure}
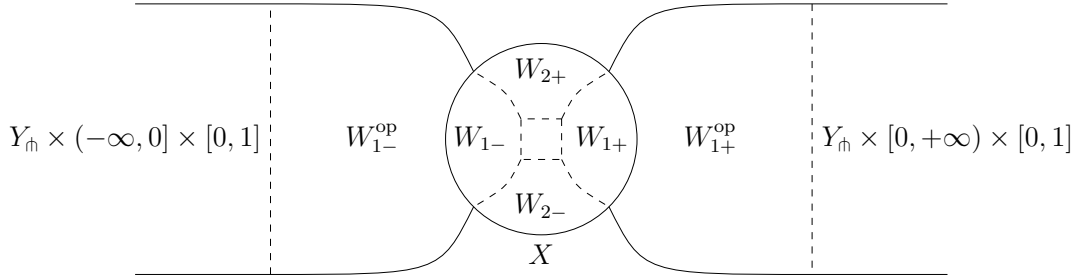
	Since the construction is similar to the one that  appeared in the proof of Theorem \ref{thm:cube}, we shall be brief. For simplicity, let us describe the construction of the metric  $\overbar g$ on $Z$ for the case where $m=2$. The general case of $m>2$ is completely similar by iterating the same construction below. 
By Sard's theorem, 	without loss of generality, we may assume that $$E_+\coloneqq\{x\in X: \varphi_1(x)=\varphi_2(x)\} \text{ and } E_-\coloneqq\{x\in X: \varphi_1(x)=-\varphi_2(x)\}$$
	 are both submanifolds in $X$.  
	For each $i = 1, 2$, the subspaces  $E_+$, $E_-$, $Y_\pitchfork\times \partial_{i\pm}I^2$  and $\partial_{i\pm}X$ together bound a  submanifold with corners in $W$, denoted by $W_{i\pm}$ (see Figure \ref{fig:W2}).  Let $W_{i\pm}^\text{op}$ be the corresponding copy of $W_{i\pm}$ in $W^\text{op}$.
	 Set $W_2$ to be the following union of subspaces in $Z$: 
	 \begin{align*}
	 W_2=&(Y_\pitchfork\times \R_{\leq 0} \times[0,1])\cup
	 W_{1-}^{\text{op}}\cup X \cup W_{1+}^\text{op}\cup (Y_\pitchfork\times\R_{\geq 0} \times[0,1])
	 \end{align*}
	 Then $W_2$ is a manifold with boundary, whose boundary consists of two components $\partial_-W_2$ and $\partial _+W_2$.
	One can extend the metric $g$ on $X$ to a metric $g_2$ on $W_2$ such that the $g_2$-distance $\dist_{g_2}(\partial_-W_2, \partial_+W_2)$ of $\partial_-W_2$ and $\partial _+W_2$ is $\geq \ell_2-\varepsilon/3$. Furthermore, we require that $\overbar g$ coincides with the product metric $g_{Y_\pitchfork}+dx_1^2 + dx_2^2$ on $Y_\pitchfork \times \mathbb R\times [0, \ell_2] $ outside a compact set of $W_2$. 
	
	We construct  $W_1$ and a metric $g_1$ on $W_1$ with similar properties.  By viewing $W_1\cup_X W_2$ as a submanifold of $Z$, we finally extend the metric on  $W_1\cup_X W_2$ to obtain a metric $\overbar g$ on $Z$.
	  
	  For $x\in Z$, let $\overbar\varphi_i(x)$ be the signed distance from $x$ to $\partial_- W_i$. Let  $z_i$   be an appropriate smooth approximation of  $\overbar\varphi_i-\ell_i/2$. To summarize, we have constructed a complete metric $\overbar g$  and real-valued smooth functions $z_1$ and $z_2$ on $Z$ with the following properties.
	  \begin{itemize}
	  	\item The metric $\overbar g$ on $Z$ restricts to the metric $g$ on $X$, and coincides with the product metric $R^2(g_{Y_\pitchfork}+dx_1^2+dx_2^2)$ outside a compact subset of $Z$ for some $R>0$. In particular, there exists a positive number $\sigma>0$ such that $\Sc_{\overbar g}\geq -\sigma$ on $Z$.
	  	\item For each $x\in Z$, if both $|z_1(x)|\leq \ell_1/2-\varepsilon$ and $|z_2(x)|\leq  \ell_2/2-\varepsilon$ for some sufficiently small $\varepsilon >0$,  then $x$ lies in $X$. 
	  	\item For each $C > 0$, the subset $\{x\in \R^2: z_1(x)^2+z_2(x)^2\leq C\}$ is compact.
	  	\item $\|\nabla z_i\|\leq 1+\varepsilon$ for both $i=1, 2$. Here $\nabla$ is the gradient with respect to the metric $\overbar g$ and $\|\nabla z_i\|$ is measured with respect to the metric $\overbar g$.
	  \end{itemize}
	  For the general case where $m\geq 2$, the same induction argument above produces a complete metric $\overbar g$ on $Z$ and real-valued  smooth functions $\{z_i\}_{1\leq i \leq n}$ on $Z$ with similar properties. By the construction of $Z$, it is not difficult to see that the $\Gamma$-covering space $\widetilde X \to X$ naturally extends to a  $\Gamma$-covering space $\widetilde Z$ of $Z$.  The metric $\overbar g$ and the functions $z_i$ on $Z$ lift to $\widetilde Z$, which will be denoted by $\widetilde g$ and $\widetilde z_i$ respectively.

	  Now let $T\widetilde Z$ be the tangent bundle of $\widetilde Z$, and $\mathbb E^m$ the trivial bundle over $\widetilde Z$. Let
	  $$E=S(T\widetilde Z\oplus \mathbb E^m)$$
	  be the spinor bundle of $T\widetilde Z\oplus \mathbb E^m$. Consider the following Callias-type operator
	  $$B=D\hotimes 1+\sum_{i=1}^m 1\hotimes f_i(\widetilde z_i)\hat e_i$$
	  acting on $L^2(\widetilde Z, E)$,  
	  where $\hat e_i$ denotes the Clifford action of $\cl_{0,m}(\mathbb E^m)$. Note that $B$ is a $\mathbb Z/2\mathbb Z$-graded $\Gamma$-equivariant essentially self-adjoint operator. 	Since we have  assumed to the contrary that
	  \begin{equation}
	  	\sum_{i=1}^{m}\frac{1}{\ell_i^2}< \frac{n^2}{4\pi^2}, 
	  \end{equation} 
	 it follows from the same estimates as in the proof of Theorem \ref{thm:cube} that $B$ is invertible. In particular, the higher index $\ind_\Gamma(B)$  of $B$ vanishes in $KO_{n-m}(C^*_{\max{}}(\Gamma; \R))$.
	  
	 	On the other hand, let $S_{\widetilde Y_\pitchfork}$ be the Clifford bundle over $\widetilde Y_\pitchfork$  with respect to the given metric $\widetilde{g}_{Y_\pitchfork}$, and $D_{\widetilde Y_\pitchfork}$ the associated Clifford-linear Dirac operator.  We consider the following generalized Bott-Dirac operator
	 	$$B_o=D_{\widetilde Y_\pitchfork}\hotimes 1\hotimes 1+\sum_{i=1}^{m}1\hotimes\frac{\partial}{\partial x_i}e_i\hotimes 1+\sum_{i=1}^{m}1\hotimes1\hotimes x_i\hat e_i $$
	 	acting on $L^2(\widetilde Y_\pitchfork\times\mathbb E^{m}, S_{\widetilde Y_\pitchfork}\hotimes S(\mathbb E^m\oplus \mathbb E^m))$. By the product formula of the higher index, we see that $$\ind_\Gamma(B_o)=\ind_\Gamma(D_{\widetilde Y_\pitchfork})\in KO_{n-m}(C^*_{\max{}}(\Gamma,\R)),$$ where $\ind_\Gamma(D_{\widetilde Y_\pitchfork})$ is assumed to be non-zero. Furthermore, by applying the relative higher index theorem \cite{Bunke, xieyu2014}, we have 
	 	$$\ind_\Gamma(B_o)-\ind_\Gamma(B)=\ind_\Gamma (D_{\mathscr D\widetilde X})$$
	 	in $KO_{n-m}(C^*_{\max{}}(\Gamma; \R))$, 
	 	where $D_{\mathscr D\widetilde X}$ is the associated Clifford-linear Dirac operator on the double  $\mathscr D \widetilde X$ of $\widetilde X$. Since $\ind_\Gamma (D_{\mathscr D\widetilde X}) = 0$ (cf.  \cite[Theorem 5.1]{xieyu2014}), it follows that 
	 	$$\ind_\Gamma(B_o)= \ind_\Gamma(B) + \ind_\Gamma (D_{\mathscr D\widetilde X}) = 0.$$  This leads to a contradiction, hence finishes the proof. 

 Now let us prove Part 2. 
 \proofpart{2}{Strict inequality}
 Recall that in the proof of Part 1 above, there is one extra step that allows us to reduce the Part 1 to the case where  $\pi_1(Y_\pitchfork) \to \pi_1(X)$ is split injective. Strictly speaking, the construction in this extra step may shrink the distances  $\ell_i = \dist(\partial_{i-}, \partial_{i+})$ by a very small amount. Such a shrinking does not affect the proof of the nonstrict inequality in Part 1. However, we are not allowed to shrink the distances  $\ell_i$ by any amount when proving the strict inequality. For this reason, we impose the slightly stronger assumption that $\pi_1(Y_\pitchfork) \to \pi_1(X)$ is split injective in Part 2. With the assumption of split injectivity, The rest of the proof proceeds similarly to that of  Part 2 of Theorem \ref{thm:cube}.

	\end{proof}


\begin{thebibliography}{10}
	
	\bibitem{MR1233861}
	N.~Anghel.
	\newblock On the index of {C}allias-type operators.
	\newblock {\em Geom. Funct. Anal.}, 3(5):431--438, 1993.
	
	\bibitem{MR3410545}
	Jean-Pierre Bourguignon, Oussama Hijazi, Jean-Louis Milhorat, Andrei Moroianu,
	and Sergiu Moroianu.
	\newblock {\em A spinorial approach to {R}iemannian and conformal geometry}.
	\newblock EMS Monographs in Mathematics. European Mathematical Society (EMS),
	Z\"{u}rich, 2015.
	
	\bibitem{Bunke}
	Ulrich Bunke.
	\newblock A {$K$}-theoretic relative index theorem and {C}allias-type {D}irac
	operators.
	\newblock {\em Math. Ann.}, 303(2):241--279, 1995.
	
	\bibitem{MR4181824}
	Simone Cecchini.
	\newblock A long neck principle for {R}iemannian spin manifolds with positive
	scalar curvature.
	\newblock {\em Geom. Funct. Anal.}, 30(5):1183--1223, 2020.
	
	\bibitem{Chodosh:2020tk}
	Otis Chodosh and Chao Li.
	\newblock Generalized soap bubbles and the topology of manifolds with positive
	scalar curvature.
	\newblock arXiv:2008.11888, 2020.
	
	\bibitem{MR532376}
	R.~E. Greene and H.~Wu.
	\newblock {$C^{\infty }$} approximations of convex, subharmonic, and
	plurisubharmonic functions.
	\newblock {\em Ann. Sci. \'{E}cole Norm. Sup. (4)}, 12(1):47--84, 1979.
	
	\bibitem{MR3822551}
	Mikhael Gromov.
	\newblock A dozen problems, questions and conjectures about positive scalar
	curvature.
	\newblock In {\em Foundations of mathematics and physics one century after
		{H}ilbert}, pages 135--158. Springer, Cham, 2018.
	
	\bibitem{MR3816521}
	Mikhael Gromov.
	\newblock Metric inequalities with scalar curvature.
	\newblock {\em Geom. Funct. Anal.}, 28(3):645--726, 2018.
	
	\bibitem{Gromov:2019aa}
	Mikhael Gromov.
	\newblock Four lectures on scalar curvature.
	\newblock arXiv:1908.10612, 2019.
	
	\bibitem{MGBL80b}
	Mikhael Gromov and H.~Blaine Lawson, Jr.
	\newblock The classification of simply connected manifolds of positive scalar
	curvature.
	\newblock {\em Ann. of Math. (2)}, 111(3):423--434, 1980.
	
	\bibitem{MGBL80}
	Mikhael Gromov and H.~Blaine Lawson, Jr.
	\newblock Spin and scalar curvature in the presence of a fundamental group.
	{I}.
	\newblock {\em Ann. of Math. (2)}, 111(2):209--230, 1980.
	
	\bibitem{MGBL83}
	Mikhael Gromov and H.~Blaine Lawson, Jr.
	\newblock Positive scalar curvature and the {D}irac operator on complete
	{R}iemannian manifolds.
	\newblock {\em Inst. Hautes {\'E}tudes Sci. Publ. Math.}, (58):83--196 (1984),
	1983.
	
	\bibitem{Gromov:2020aa}
	Misha Gromov.
	\newblock No metrics with positive scalar curvatures on aspherical 5-manifolds.
	\newblock 2020.
	
	\bibitem{GXY2020}
	Hao Guo, Zhizhang Xie, and Guoliang Yu.
	\newblock Quantitative {K}-theory, positive scalar curvature, and band width.
	\newblock {\em to appear in special volume Perspectives on Scalar Curvature
		{(editors: Gromov and Lawson)}}, 2020.
	\newblock arXiv:2010.01749.
	
	\bibitem{Lohkamp:2018wp}
	Joachim Lohkamp.
	\newblock Minimal smoothings of area minimizing cones.
	\newblock  2018.
	
	\bibitem{RSSY79b}
	R.~Schoen and S.~T. Yau.
	\newblock On the structure of manifolds with positive scalar curvature.
	\newblock {\em Manuscripta Math.}, 28(1-3):159--183, 1979.
	
	\bibitem{RSSY79}
	R.~Schoen and Shing~Tung Yau.
	\newblock Existence of incompressible minimal surfaces and the topology of
	three-dimensional manifolds with nonnegative scalar curvature.
	\newblock {\em Ann. of Math. (2)}, 110(1):127--142, 1979.
	
	\bibitem{Schoen:2017aa}
	Richard Schoen and Shing-Tung Yau.
	\newblock Positive scalar curvature and minimal hypersurface singularities.
	\newblock arXiv 1704.05490, 2017.
	
	\bibitem{MR1629775}
	Wolfgang Walter.
	\newblock {\em Ordinary differential equations}, volume 182 of {\em Graduate
		Texts in Mathematics}.
	\newblock Springer-Verlag, New York, 1998.
	\newblock Translated from the sixth German (1996) edition by Russell Thompson,
	Readings in Mathematics.
	
	\bibitem{WangXie2022}
	Jinmin Wang and Zhizhang Xie.
	\newblock Dihedral rigidity for submanifolds of warped product manifolds.
	\newblock arXiv:2303.13492.
	
	\bibitem{Wang:2021uj}
	Jinmin Wang, Zhizhang Xie, and Guoliang Yu.
	\newblock Decay of scalar curvature on uniformly contractible manifolds with
	finite asymptotic dimension.
	\newblock {\em \emph{to appear in} Communications on Pure and Applied
		Mathematics}, 2021.
	\newblock arXiv:2101.11584.
	
	\bibitem{Wang:2021tq}
	Jinmin Wang, Zhizhang Xie, and Guoliang Yu.
	\newblock On {G}romov's dihedral extremality and rigidity conjectures.
	\newblock arXiv:2112.01510, 2021.
	
	\bibitem{MR1501735}
	Hassler Whitney.
	\newblock Analytic extensions of differentiable functions defined in closed
	sets.
	\newblock {\em Trans. Amer. Math. Soc.}, 36(1):63--89, 1934.
	
	\bibitem{Xie:2021tm}
	Zhizhang Xie.
	\newblock A quantitative relative index theorem and {G}romov's conjectures on
	positive scalar curvature.
	\newblock 2021.
	\newblock to appear in \emph{Journal of Noncommutative Geometry}.
	
	\bibitem{xieyu2014}
	Zhizhang Xie and Guoliang Yu.
	\newblock A relative higher index theorem, diffeomorphisms and positive scalar
	curvature.
	\newblock {\em Adv. Math.}, 250:35--73, 2014.
	
	\bibitem{Zeidler:2019aa}
	Rudolf Zeidler.
	\newblock Band width estimates via the {D}irac operator.
	\newblock {\em Journal of Differential Geometry}, 2020.
	
	\bibitem{MR4181525}
	Rudolf Zeidler.
	\newblock Width, largeness and index theory.
	\newblock {\em SIGMA Symmetry Integrability Geom. Methods Appl.}, 16:Paper No.
	127, 15, 2020.
	
\end{thebibliography}
\end{document}